%% file: main.tex
\documentclass[hidelinks,onefignum,onetabnum]{siamart220329}

\usepackage{cancel}
\input{art_shared}

\ifpdf
\hypersetup{
  pdftitle={AP and energy stable dynamical low-rank approximation},
  pdfauthor={}
}
\fi

\usepackage{bm}
\usepackage{subcaption}

\newcommand{\bg}{\mathbf{g}}

\newcommand{\ba}{\mathbf{a}}

\usepackage{calligra,amsmath,amssymb}

\begin{document}

\maketitle

\begin{abstract}
Radiation transport problems are posed in a high-dimensional phase space, limiting the use of finely resolved numerical simulations. An emerging tool to efficiently reduce computational costs and memory footprint in such settings is dynamical low-rank approximation (DLRA). Despite its efficiency, numerical methods for DLRA need to be carefully constructed to guarantee stability while preserving crucial properties of the original problem. Important physical effects that one likes to preserve with DLRA include capturing the diffusion limit in the high-scattering regimes as well as dissipating energy. In this work we propose and analyze a dynamical low-rank method based on the ``unconventional" basis-update \& Galerkin step integrator. We show that this method is asymptotic--preserving, i.e., it captures the diffusion limit, and energy stable under a CFL condition. The derived CFL condition captures the transition from the hyperbolic to the parabolic regime when approaching the diffusion limit. 
\end{abstract}



\section{Introduction}
Radiative transfer equations describe the dynamics of radiation traveling through a background medium. The interplay of particle advection as well as collisions with the background material yields the integro-differential radiation transport equation for the particle density (also called angular flux) $f$. Assuming mono-energetic particles, this particle density depends on time, space and direction of travel, i.e., for a three-dimensional spatial domain, the phase-space is at least six-dimensional, resulting in prohibitive memory and computational costs when a fine discretization is chosen. Conventional numerical methods are therefore frequently constructed to approximate the angular flux with a sufficiently coarse resolution while accounting for arising numerical artifacts through, for example, the use of filters or ray-effect mitigation. Moreover, numerical methods are often designed to capture the diffusion limit. When scattering is large on a small time scale, the radiative transfer equation tends towards a parabolic diffusion equation. Methods conserving this behavior and efficiently treat arising stiffness when scattering terms become large are called asymptotic--preserving (AP) \cite{jin1999efficient} methods. There has been significant development of AP methods for hyperbolic and kinetic equations in the past twenty years. The readers are referred to \cite{Jin_Rev,HJL17,Jin} for thorough review of these methods.


To tackle computational challenges and reduce memory, dynamical low-rank approximation \cite{koch2007dynamical} has been employed for various problems in radiative transfer including for example radiation therapy \cite{kusch2021robust}, high-scattering regimes \cite{ding2019error,einkemmer2021asymptotic,ding2021dynamical} and criticality of nuclear systems \cite{kusch2022low}. Further works where DLRA has proven to yield an efficient method in the context of radiative transfer are \cite{peng2020low,peng2021high}. The core idea of dynamical low-rank approximation is to constrain the evolution of the solution to a rank $r$ manifold $\mathcal{M}$ by projecting the dynamics onto the tangent plane of $\mathcal{M}$. Intuitively, the resulting evolution equations can be interpreted as a Galerkin system with $r$ basis functions for each phase-space dimension, which updates coefficients and basis functions in time according to the dynamics of the problem. Robust integrators for the corresponding evolution equations have been proposed in \cite{lubich2014projector,ceruti2022unconventional}. The projector--splitting integrator (PSI) \cite{lubich2014projector} splits the projection onto the tangent plane into three subflows which are then solved consecutively. Though exhibiting robustness irrespective of small singular values within the solution \cite{kieri2016discretized}, the PSI includes a backward in time step, which can yield unstable evolution equations in the case of parabolic and even hyperbolic problems \cite{kusch2021stability}. Moreover, as we show in this work, propagating the solution backward in time can lead to an unphysical increase of the solution's energy. 

The ``unconventional" basis update \& Galerkin step (BUG) integrator \cite{ceruti2022unconventional} first updates the basis functions and then performs a Galerkin step to evolve the expansion coefficients. While sharing the robustness property of the PSI \cite{ceruti2022unconventional}, it does not require a step backward in time and allows the basis functions to be updated in parallel. Moreover, it allows for an efficient way to augment the basis \cite{ceruti2022rank}, which has been employed to propose a robust rank adaptive BUG integrator \cite{ceruti2022rank} that has been extended to tree tensor networks in \cite{ceruti2022ranktree}. Moreover, the basis augmentation allows for a construction of a conservative extension to the BUG integrator \cite{einkemmer2022robust}.

In this work, we propose an asymptotic--preserving BUG integrator. Similar to the asymptotic--preserving DLRA scheme proposed in \cite{einkemmer2021asymptotic}, our scheme is constructed by a micro-macro decomposition \cite{lemou2008new}. In contrast to \cite{einkemmer2021asymptotic}, choosing the BUG integrator to evolve the low-rank factors, we are able to prove stability under a CFL condition, which captures the asymptotic behavior. I.e., for Knudsen numbers of $O(1)$, we obtain a hyperbolic CFL condition whereas for small Knudsen numbers the time step restriction becomes a parabolic condition. Moreover, by choosing a modal discretization in angle, we are able to obtain a symmetric flux matrix, which ensures hyperbolicity. The following novel results demonstrate the beneficial properties of the proposed scheme:
\begin{itemize}
\item \textit{An accurate stability analysis for the BUG integrator}: We present an analysis for the micro-macro P$_N$ equations that captures the diffusive limit. While such an analysis has been proposed for nodal discretizations \cite{lemou2008new}, the extension to modal schemes is not straightforward and to the best of the authors' knowledge has not been shown before. Note that for the micro-macro decomposition the P$_N$ and S$_N$ methods are not equivalent, which makes the diagonalization of the scheme difficult. We enable the diagonalization of the scheme by proposing a modified stabilization. 
\item \textit{A stability analysis for the proposed dynamical low-rank scheme}: We show that the BUG integrator allows for a stability analysis, which is based on the full micro--macro equations.
\item \textit{A derivation of the AP property for the BUG integrator}: Based on \cite{einkemmer2021asymptotic}, we show that the proposed scheme is asymptotic--preserving.
\end{itemize}
The paper is structured as follows: After the introduction, we provide the necessary background on computational methods for radiative transfer in Section~\ref{sec:recapRadTraf} and dynamical low-rank approximation in Section~\ref{sec:recapDLRA}. The modal micro-macro scheme in proposed in Section~\ref{sec:discFull} along with a proof of energy stability. In Section~\ref{sec:DLRAmicromacro} we extend the proposed scheme to dynamical low-rank approximation and show that the scheme is asymptotic--preserving. The discretization for the proposed DLRA equations is derived in Section~\ref{sec:DLRAdisc} and is shown to preserve the correct energy dissipation rate in a time-continuous setting. Moreover, we show energy stability in the fully discrete setting. Lastly, we present numerical results in Section~\ref{sec:numres}.
\newpage
\section{Background}

\subsection{Recap: Radiative transfer}\label{sec:recapRadTraf}
In the following, we briefly review the multiscale radiative transfer equation, the micro-macro decomposition \cite{lemou2008new} as well as its spherical harmonics (P$_N$) approximation \cite{case1967linear}. For simplicity of presentation, we focus on the one-dimensional radiative transfer equation in diffusive scaling
\begin{equation} \label{eq:radTransport}
\partial_t f+\frac{1}{\varepsilon} \mu \partial_x f=\frac{\sigma(x)}{\varepsilon^2}\left(\frac1{\sqrt{2}}\rho -f\right).
\end{equation}
The scalar flux $f:\mathbb{R}_+\times D \times [-1,1] \rightarrow \mathbb{R}_+$ depends on time $t$, space $x\in D\subset \mathbb{R}$ and direction of travel $\mu\in [-1,1]$. The above equation is equipped with initial and boundary conditions which we will state for the individual problems later in this work. The scattering cross-section $\sigma : D\rightarrow \mathbb{R}_+$ models the probability of particles scattering with the background material, where $\sigma(x)\geq \sigma_0 > 0$. The scalar flux is given by $\rho := \frac1{\sqrt{2}}\langle f\rangle_{\mu}$, where the integration over the angular domain is denoted by $\langle \cdot \, \rangle_{\mu}$. Letting the Knudsen number $\varepsilon$ go to zero increases the number of scattering events per unit time. In the limit $\varepsilon\rightarrow 0$ the angular flux $f$ tends towards the scalar flux $\rho$ and the radiative transfer equation changes to a parabolic diffusion equation of the form 
\begin{align*}
\partial_t \rho - \frac13\partial_x \left( \frac{1}{\sigma(x)}\partial_x \rho(x)  \right) = 0.
\end{align*}
A frequently employed method to discretize the directional variable $\mu$ is the P$_N$ method \cite{case1967linear}. In our one-dimensional setting, the P$_N$ ansatz reads
\begin{equation}
f(t,x,\mu)\approx f_N(t,x,\mu):=\sum_{l=0}^N f_l(t,x)p_l(\mu),
\end{equation}
where $\{p_l\}_{l=0}^N$ are the orthonormal Legendre polynomials. Since $\langle p_l p_k\rangle=\delta_{lk}$ we readily obtain $f_l=\langle fp_l\rangle$. 
Note that $p_l$ satisfies the recurrence relation
\begin{equation}
\mu p_l(\mu)=a_{l-1}p_{l-1}(\mu)+a_lp_{l+1}(\mu), \quad a_l=\frac{l+1}{\sqrt{(2l+1)(2l+3)}}.
\end{equation}
Substituting $f_N$ into \eqref{eq:radTransport} and projecting $\langle \cdot \ p_k\rangle$, one obtains the P$_N$ moment equations
\begin{equation}
\left\{
\begin{split}
&\partial_t f_0+\frac{1}{\varepsilon}a_0\partial_xf_1=0,\\
&\partial_t f_l+\frac{1}{\varepsilon}(a_{l-1}\partial_xf_{l-1}+a_l\partial_x f_{l+1})=-\frac{\sigma}{\varepsilon^2}f_l, \quad 1\leq l\leq N-1,\\
&\partial_tf_N+\frac{1}{\varepsilon}a_{N-1}\partial_x f_{N-1}=-\frac{\sigma}{\varepsilon^2}f_N.
\end{split}
\right.
\end{equation}
To facilitate the construction of asymptotic--preserving schemes, the micro-macro decomposition has been proposed in \cite{lemou2008new}. 
The key ingredient of this decomposition is to write the solution as the scalar flux $\rho:=f_0=\frac{1}{\sqrt{2}}\langle f\rangle$ plus the microscopic correction term
\begin{equation}
\bg:=\frac{1}{\varepsilon}(f_1,\dots,f_N)^{\top}.
\end{equation}
Hence, with $\mathbf{f} = \rho\mathbf{e}_1 + \varepsilon \,(0, g_1,\cdots,g_N)^{\top}$, where $\mathbf{e}_1\in\mathbb{R}^{N+1}$ is the first Euclidean unit vector $(e_1)_i = \delta_{i1}$, the above PN system can be written as
\begin{equation} \label{eq:PN}
\left\{
\begin{split}
&\partial_t \rho+a_0\partial_x g_1=0,\\
&\partial_t \bg+\frac{1}{\varepsilon}\mathbf A\partial_x \bg=-\frac{1}{\varepsilon^2}\ba \partial_x \rho-\frac{\sigma}{\varepsilon^2}\bg,
\end{split}
\right.
\end{equation}
where
\begin{equation}
\mathbf A=
\begin{pmatrix}
0& a_1 &  &  \\
a_1 & \ddots &\ddots &  \\
  & \ddots  & \ddots & a_{N-1}  \\
 &  & a_{N-1} & 0
\end{pmatrix}_{N\times N},
\quad  \ba =
\begin{pmatrix}
a_0\\
0\\
\vdots\\
0
\end{pmatrix}_{N\times 1}.
\end{equation}

\subsection{Recap: Dynamical low-rank approximation}\label{sec:recapDLRA}
In this section, we provide a brief review of the dynamical low-rank approximation proposed in \cite{koch2007dynamical}. The core idea of a dynamical low-rank approximation (DLRA) is to evolve the solution on a low-rank manifold, e.g., given a spatially discretized solution $g_{lk} := g_l(t,x_k)$ we wish to evolve $\mathbf g$ such that it remains on the manifold of rank $r$ matrices $\mathcal{M}_r$. Hence, every solution $\mathbf{g}_r\in \mathcal{M}_r\subset\mathbb{R}^{N_x\times N}$ can be written as 
\begin{align}\label{eq:rankrsol}
\mathbf{g}_r = \mathbf{X}(t)\mathbf{S}(t)\mathbf{V}(t)^{\top}.
\end{align}
That is, the solution is spanned by the spatial basis matrix $\mathbf{X}:\mathbb{R}_+\rightarrow \mathbb{R}^{N_x \times r}$ and the moment basis matrix $\mathbf{V}:\mathbb{R}_+\rightarrow \mathbb{R}^{N \times r}$. The coefficient matrix is given by $\mathbf{S}:\mathbb{R}_+\rightarrow \mathbb{R}^{r \times r}$. To evolve basis matrices and coefficients in time, evolution equations are derived such that for $\mathbf g_r\in\mathcal{M}_r$ and a given right-hand side $\mathbf{F}$ we impose
\begin{align}\label{eq:DLRproblem}
\dot{\mathbf g}_r(t)\in T_{ \mathbf g_r(t)}\mathcal{M}_r \qquad \text{such that} \qquad \left\Vert \dot{\mathbf g}_r(t)-\mathbf{F}(\mathbf g_r(t)) \right\Vert = \text{min}.
\end{align}
The tangent space of $\mathcal{M}_r$ at $\mathbf g_r(t)$ is denoted by $T_{\mathbf g_r(t)}\mathcal{M}_r$ and the norm $\Vert\cdot\Vert$ denotes the Frobenius norm. Condition \eqref{eq:DLRproblem} can be reformulated \cite[Lemma~4.1]{koch2007dynamical} as 
\begin{align}\label{eq:lowRankProjector}
\dot{\mathbf g}_r(t) = \mathbf{P}(\mathbf g_r(t))\mathbf{F}(\mathbf g_r(t)),
\end{align}
where $\mathbf{P}$ is the orthogonal projection onto the tangent space
\begin{align*}
\mathbf P\mathbf h = \mathbf{X}\mathbf{X}^{\top} \mathbf h - \mathbf{X}\mathbf{X}^{\top} \mathbf h \mathbf{V}\mathbf{V}^{\top} + \mathbf h \mathbf{V}\mathbf{V}^{\top}.
\end{align*}
Following \cite{koch2007dynamical}, evolution equations of the factorized solution can be derived from the above equation as follows
\begin{subequations}\label{eq:equationsFactors}
\begin{align}
    \dot{\mathbf{S}}(t) =& \mathbf{X}(t)^\top\mathbf{F}(\mathbf g_r(t))\mathbf{V}(t), \\
    \dot{\mathbf{X}}(t) =& (\mathbf I - \mathbf{X}(t)\mathbf{X}(t)^\top)\mathbf{F}(\mathbf g_r(t))\mathbf{V}(t)\mathbf{S}(t)^{-1}, \\
    \dot{\mathbf{V}}(t) =& (\mathbf I - \mathbf{V}(t)\mathbf{V}(t)^\top)\mathbf{F}(\mathbf g_r(t))^\top\mathbf{X}(t)\mathbf{S}(t)^{-\top}.
\end{align}
\end{subequations}
To avoid the inversion of the coefficient matrix on the right-hand side, robust integrators have been developed \cite{lubich2014projector,ceruti2022unconventional}. 

While the matrix projector--splitting integrator \cite{lubich2014projector} includes a step backward in time, the ``unconventional'' basis update \& Galerkin step (BUG) integrator \cite{ceruti2022unconventional} evolves the solution forward in every step. The BUG integrator evolves a given factorized solution at time $t_0$ to time $t_1$ according to the following scheme: 
\begin{enumerate}
    \item \textbf{$K$-step}: Update $\mathbf X^{0}$ to $\mathbf X^{1}$ via
    \begin{align}
        \dot{\mathbf K}(t) &= \mathbf{F}(\mathbf{K}(t)\mathbf{V}^{0,T})\mathbf{V}^0, \qquad \mathbf K(t_0) = \mathbf{X}^0\mathbf{S}^0.\label{eq:KStepSemiDiscreteUI}
    \end{align}
Perform a QR decomposition to obtain $\mathbf K(t_1) = \mathbf X^1 \mathbf R$ and $\mathbf M = \mathbf X^{1,T}\mathbf X^0$. Note that $\mathbf R$ is thrown away after this step.
\item \textbf{$L$-step}: Update $\mathbf V^0$ to $\mathbf V^1$ via
\begin{align}
\dot{\mathbf L}(t) &= \mathbf{X}^{0,T}\mathbf{F}(\mathbf{X}^0\mathbf{L}(t)), \qquad \mathbf L(t_0) = \mathbf{S}^0 \mathbf{V}^{0,T}.\label{eq:LStepSemiDiscreteUI}
\end{align}
Perform a QR decomposition to obtain $\mathbf L(t_1) = \mathbf V^1\mathbf{\widetilde R}$ and $\mathbf N = \mathbf V^{1,T} \mathbf V^0$. Note that $\mathbf{\widetilde R}$ is thrown away after this step.
\item \textbf{$S$-step}: Update $\mathbf S^0$ to $\mathbf S^1$ via
\begin{align}
\dot{\mathbf S}(t) = \mathbf{X}^{1,T}\mathbf{F}(\mathbf{X}^{1}\mathbf{S}(t)\mathbf{V}^{1,T})\mathbf{V}^1, \qquad \mathbf S(t_0) &= \mathbf M\mathbf S^0 \mathbf N^{\top}\label{eq:SStepSemiDiscreteUI}
\end{align}
and set $\mathbf S^1 = \mathbf S(t_1)$.
\end{enumerate}
The time updated solution is then given by $\mathbf{g}(t_1) = \mathbf{X}^1\mathbf{S}^1\mathbf{V}^{1,T}$. Note that the BUG integrator, while being able to compute $K$ and $L$-steps in parallel, is only first order accurate in time.

\section{Energy stable and modal micro-macro discretization}\label{sec:discFull}
Before deriving evolution equations for the factorized solution of the BUG integrator, we present an asymptotic--preserving discretization to the full micro-macro P$_N$ system \eqref{eq:PN}.
To derive a suitable stabilization matrix, note that the flux matrix $\mathbf A$ can be written in terms of Gauss-Legendre quadrature rule as
\begin{align}
    a_{ij} = \langle \mu p_i p_j \rangle = \sum_{k=1}^{N+1}w_k p_i(\mu_k)p_j(\mu_k)\mu_k.
\end{align}
That is, defining the transformation matrix $\mathbf{T}\in\mathbb{R}^{N\times (N+1)}$ with entries $T_{ik} = \sqrt{w_k}p_i(\mu_k)$ and $\mathbf{M} = \text{diag}(\mu_1,\cdots,\mu_{N+1})$ we have $\mathbf{A} = \mathbf{T}\mathbf{M}\mathbf{T}^{\top}$. With this definition, we define the stabilization matrix $|\mathbf{A}| := \mathbf{T}|\mathbf{M}|\mathbf{T}^{\top}$ as well as $\mathbf{A}^{\pm} := \frac12\mathbf{T}(\mathbf{M} \pm|\mathbf{M}|)\mathbf{T}^{\top}$. 
\begin{remark}\label{rem:stabilization}
Note that this choice of the stabilization matrix is not equivalent to the commonly used Roe matrix $\mathbf{\widetilde T}|\mathbf{\widetilde M}|\mathbf{\widetilde T}^{\top}$ where $\mathbf{A} = \mathbf{\widetilde T}\mathbf{\widetilde M}\mathbf{\widetilde T}^{\top}$ is the eigendecomposition. That is, $\mathbf{\widetilde T}\in\mathbb{R}^{N\times N}$ collects the orthonormal eigenvectors of $\mathbf{A}$ and $\mathbf{\widetilde M}\in\mathbb{R}^{N\times N}$ is the diagonal eigenvalue matrix. Instead, the derived factorization of the flux matrix uses transformation matrices that are elements of $\mathbb{R}^{N\times (N+1)}$, which is needed to later diagonalize the scheme when proving energy stability. 
\end{remark}
For the spatial discretization, we define an equidistant grid with cell interface points $x_{1/2} \leq \cdots \leq x_{N_x + 1/2}$ and midpoints $x_j$ where $j\in\{1,\cdots,N_x\}$. Then, the scalar flux is represented on the midpoint $x_j$ at time $t_n$ as $\rho_j^{n}$. The microscopic correction is presented on the cell interfaces at time $t_n$ as $\bg_{j+1/2}^{n}$.
Then, the fully discretized system \eqref{eq:PN} reads as
\begin{equation}\label{eq:micromacro_orig}
\left\{
\begin{split}
&\frac{\rho_j^{n+1}-\rho_j^n}{\Delta t}+a_0\frac{g_{1,j+1/2}^{n+1}-g_{1,j-1/2}^{n+1}}{\Delta x}=0,\\
&\frac{\bg_{j+1/2}^{n+1}-\bg_{j+1/2}^n}{\Delta t}+\frac{1}{\varepsilon}\frac{\mathbf A^+(\bg_{j+1/2}^n-\bg_{j-1/2}^n)+\mathbf A^-(\bg_{j+3/2}^n-\bg_{j+1/2}^n)}{\Delta x},\\
 &\qquad =-\frac{1}{\varepsilon^2}\ba \frac{\rho_{j+1}^n-\rho_j^n}{\Delta x}-\frac{\sigma_{j+1/2}}{\varepsilon^2}\bg^{n+1}_{j+1/2}.
\end{split}
\right.
\end{equation}
To arrive at a more compact notation, we define the discretized advection operator 
\begin{align}
    \mathcal{L}\bg_{j+1/2}^n := (\mathbf A^+\mathcal{D}^-+\mathbf A^-\mathcal{D}^{+})\bg_{j+1/2}^n
\end{align}
with the upwind stencil terms
\begin{equation}
\mathcal{D}^-\bg_{j+1/2}=\frac{\bg_{j+1/2}-\bg_{j-1/2}}{\Delta x}, \quad \mathcal{D}^{+}\bg_{j+1/2}=\frac{\bg_{j+3/2}-\bg_{j+1/2}}{\Delta x}.
\end{equation}
Then, the discrete micro-macro system \eqref{eq:micromacro_orig} reads
\begin{equation} \label{eq:schemeFull}
\left\{
\begin{split} 
&\frac{\rho_j^{n+1}-\rho_j^n}{\Delta t}+a_0\mathcal{D}^-g_{1,j+1/2}^{n+1}=0,\\
&\frac{\bg_{j+1/2}^{n+1}-\bg_{j+1/2}^n}{\Delta t}+\frac{1}{\varepsilon}\mathcal{L}\bg_{j+1/2}^n=-\frac{1}{\varepsilon^2}\ba \mathcal{D}^{+}\rho_{j}^n-\frac{\sigma_{j+1/2}}{\varepsilon^2}\bg^{n+1}_{j+1/2}.
\end{split}
\right.
\end{equation}
To investigate stability, we define the energy of the solution as 
\begin{align*}
    e^n := \|\rho^n\|^2 + \varepsilon^2  \|\bg^n\|^2,
\end{align*}
where we define the discrete L$^2$ norm of $\rho$ and $\bg$ as
\begin{equation}
\|\rho\|^2=\sum_j (\rho_j)^2\Delta x, \quad \|\bg\|^2=\sum_j (\bg_{j+1/2}^{\top}\bg_{j+1/2}) \Delta x.
\end{equation}
It can be shown that the discrete micro-macro system \eqref{eq:schemeFull} dissipates the energy
with, depending on the regime, either a hyperbolic or parabolic time-step restriction:
\begin{theorem}[Energy stability]\label{th:estabilityFull} 
Assume that the time step size $\Delta t$ fulfills the CFL condition 
\begin{align}\label{eq:CFLcondition}
    \Delta t \leq C\left(\varepsilon\Delta x + \sigma_0\Delta x^2\right).
\end{align}
Then, the scheme \eqref{eq:schemeFull} is energy stable, that is, $e^{n+1}\leq e^n$.
\end{theorem}

To prove this Theorem, we first note several properties of the chosen discretization.
\begin{lemma}[Summation by parts]\label{le:sumbyparts}
For vectors $\bg_{j+1/2}, \mathbf{k}_{j+1/2}\in\mathbb{R}^N$ where $j = 1,\cdots,N_x$, the equality
\begin{equation}
\sum_j \mathbf{k}^{\top}_{j+1/2}D^{\pm} \bg_{j+1/2}=-\sum_j (D^{\mp}\mathbf{k}_{j+1/2})^{\top}\bg_{j+1/2}
\end{equation}
holds.
\end{lemma}
\begin{proof}
The result directly follows from the definition of $D^{\pm}$ and an index shift in the sum over spatial cells. For $\mathcal{D}^{+}$ we have
    \begin{align*}
\sum_j \mathbf{k}^{\top}_{j+1/2}\frac{\bg_{j+3/2}-\bg_{j+1/2}}{\Delta x}=&\frac1{\Delta x}\sum_j \mathbf{k}^{\top}_{j+1/2}\bg_{j+3/2}-\frac1{\Delta x}\sum_j \mathbf{k}^{\top}_{j+1/2}\bg_{j+1/2}\\
=& \frac1{\Delta x}\sum_j \mathbf{k}^{\top}_{j-1/2}\bg_{j+1/2}-\frac1{\Delta x}\sum_j \mathbf{k}^{\top}_{j+1/2}\bg_{j+1/2}\\
=& -\sum_j (\mathcal{D}^{-}\mathbf{k}_{j+1/2})^{\top}\bg_{j+1/2}.
\end{align*}
\end{proof}

In the following, we relate the micro-macro system to the full P$_N$ system with flux matrix $\mathbf{A}_f = \left(\langle p_{i-1} p_{j-1} \mu \rangle\right)_{i,j=1}^{N+1}$ and Roe matrix $|\mathbf{A}_f| = \mathbf{T}_f|\mathbf{M}|\mathbf{T}_f^{\top}$. Here, $\mathbf{T}_f = \left(\sqrt{w_k}p_{i-1}(\mu_k)\right)_{i,j=1}^{N+1}$ such that $\mathbf{A}_f = \mathbf{T}_f\mathbf{M}\mathbf{T}_f^{\top}$. Moreover, define 
\begin{align}
    \mathbf{a}_f = (0,a_0,0,\cdots,0)^{\top}.
\end{align}
Then, the following lemma holds.
\begin{lemma}[P$_N$ preservation]\label{le:AtoAf}
For a given vector $\mathbf{g}\in\mathbb{R}^N$ define its extension $\mathbf{h} := (0,g_1,\cdots,g_N)^{\top}\in\mathbb{R}^{N+1}$ as well as $\mathbf{\widehat h}_{j+1/2} := \mathbf{T}_f^{\top}\mathbf{h}_{j+1/2}\in\mathbb{R}^{N+1}$. Then,
\begin{align*}
    \mathbf{g}^{\top} \mathbf{A}^2 \mathbf g = \mathbf{\widehat h}^{\top} \mathbf M^2 \mathbf{\widehat h}, \quad \mathbf{g}^{\top} |\mathbf{A}| \mathbf g = \mathbf{\widehat h}^{\top} |\mathbf{M}| \mathbf{\widehat h}, \quad\mathbf{g}^{\top} \ba\ba^{\top} \mathbf g = \mathbf{\widehat h}^{\top} \mathbf{T}_f^{\top}\ba_f\ba_f^{\top}\mathbf{T}_f \mathbf{\widehat h}.
\end{align*}
\end{lemma}
\begin{proof}
We directly have that $\mathbf{h}^{\top}\mathbf{A}_f^2\mathbf{h} = \mathbf{\widehat h}^{\top} \mathbf M^2 \mathbf{\widehat h}$ and $\mathbf{h}^{\top}|\mathbf{A}_f|\mathbf{h} = \mathbf{\widehat h}^{\top} |\mathbf M| \mathbf{\widehat h}$. Moreover, we have
\begin{align*}
    \mathbf{A}_f\mathbf{h} = \begin{pmatrix}
        \mathbf{a}_f^{\top}\mathbf{h} \\
        \mathbf{A}\mathbf{g}
    \end{pmatrix},\quad |\mathbf{A}_f|\mathbf{h} = \begin{pmatrix}
        (|\mathbf{A}_f|\mathbf e_1)^{\top}\mathbf{h} \\
        |\mathbf{A}|\mathbf{g}
    \end{pmatrix}, \quad \mathbf{a}_f^{\top}\mathbf h = \mathbf a^{\top}\mathbf{g},
\end{align*}
where the second equality holds due to the choice of the stabilization matrix (cf. Remark~\ref{rem:stabilization}).
Hence, we have
\begin{align*}
    &\left(\mathbf{A}_f\mathbf{h}\right)^{\top}\left(\mathbf{A}_f\mathbf{h}\right) = \mathbf{h}^{\top}\mathbf{a}\mathbf{a}^{\top}\mathbf{h} + \mathbf{g}^{\top}\mathbf{A}^2\mathbf{g} = \mathbf{g}^{\top}\mathbf{A}^2\mathbf{g},\\
    &\mathbf{h}^{\top}|\mathbf{A}_f|\mathbf{h} = h_0 (|\mathbf{A}_f|\mathbf e_1)^{\top}\mathbf{h} + \mathbf{g}^{\top}|\mathbf{A}|\mathbf{g} = \mathbf{g}^{\top}|\mathbf{A}|\mathbf{g},\\
    &\mathbf h^{\top}\mathbf{a}_f\mathbf{a}_f^{\top}\mathbf h = \mathbf g^{\top}\mathbf{a}\mathbf{a}^{\top}\mathbf g= \mathbf{\widehat h}^{\top} \mathbf{T}_f^{\top}\ba_f\ba_f^{\top}\mathbf{T}_f \mathbf{\widehat h}.
\end{align*}
\end{proof}

When proving energy stability, we make use of two main properties of the advection operator $\mathcal{L}$:
\begin{lemma}[Positivity]\label{le:3.2}
For a given discrete function $\bg_{j+1/2}^n$, the advection operator fulfills the properties
\begin{align}\label{eq:le3a}
    \sum_j \bg_{j+1/2}^{n+1,\top} \mathcal{L}\bg_{j+1/2}^{n+1} =& \sum_j \frac{\Delta x}{2}\mathcal{D}^{+}\bg_{j+1/2}^{n+1,\top}|\mathbf A|\mathcal{D}^{+}\bg_{j+1/2}^{n+1}\geq 0
\end{align}
and
\begin{align}\label{eq:le3b}
    \sum_j \bg_{j+1/2}^{n+1,\top} \mathcal{L}\bg_{j+1/2}^n =& \sum_j \frac{\Delta x}{2}\mathcal{D}^{+}\bg_{j+1/2}^{n+1,\top}|\mathbf A|\mathcal{D}^{+}\bg_{j+1/2}^{n+1} \nonumber \\
    &+\sum_j(\bg_{j+1/2}^{n+1}-\bg_{j+1/2}^n)^{\top} \mathcal{L} \bg_{j+1/2}^{n+1}.
\end{align}
\end{lemma}
\begin{proof}
We first show \eqref{eq:le3a}. Let us start by noting that we can rewrite the advection operator as
\begin{align*}
    \mathcal{L}\bg_{j+1/2}=\mathbf A\frac{\bg_{j+3/2}-\bg_{j-1/2}}{2\Delta x}-\frac{\Delta x|\mathbf A|}{2}\mathcal{D}^{-}\mathcal{D}^{+}\bg_{j+1/2}.
\end{align*}
Then, we have that
\begin{align*}
    \sum_j \bg_{j+1/2}^{n+1,\top} \mathcal{L}\bg_{j+1/2}^{n+1} =&\sum_j \bg_{j+1/2}^{n+1,\top}\mathbf A\frac{\bg_{j+3/2}^{n+1}-\bg_{j-1/2}^{n+1}}{2\Delta x}-\sum_j \frac{\Delta x}{2}\bg_{j+1/2}^{n+1,\top}|\mathbf A|\mathcal{D}^{-}\mathcal{D}^{+}\bg_{j+1/2}^{n+1}.
\end{align*}
The first term on the right-hand side is zero since (due to a shift of index and symmetry of $\mathbf{A}$) 
\begin{align*}
    \sum_j \bg_{j+1/2}^{n+1,\top}\mathbf A\left(\bg_{j+3/2}^{n+1}-\bg_{j-1/2}^{n+1}\right) = \sum_j \bg_{j+1/2}^{n+1,\top}\mathbf A\bg_{j+3/2}^{n+1} - \sum_j \bg_{j+3/2}^{n+1,\top}\mathbf A\bg_{j+1/2}^{n+1} = 0.
\end{align*}
According to summation by parts, Lemma~\ref{le:sumbyparts}, the second term can be rewritten as 
\begin{align*}
    -\sum_j \frac{\Delta x}{2}\bg_{j+1/2}^{n+1,\top}|\mathbf A|\mathcal{D}^{-}\mathcal{D}^{+}\bg_{j+1/2}^{n+1} = \sum_j \frac{\Delta x}{2}\mathcal{D}^{+}\bg_{j+1/2}^{n+1,\top}|\mathbf A|\mathcal{D}^{+}\bg_{j+1/2}^{n+1}.
\end{align*}
Hence, \eqref{eq:le3a} holds. To show \eqref{eq:le3b}, we write
\begin{align*}
    \sum_j \bg_{j+1/2}^{n+1,\top} \mathcal{L}\bg_{j+1/2}^n
    =&\sum_j \bg_{j+1/2}^{n+1,\top} \mathcal{L}\bg_{j+1/2}^{n+1}-\sum_j \bg_{j+1/2}^{n+1,\top} \mathcal{L}(\bg_{j+1/2}^{n+1}-\bg_{j+1/2}^n)
\end{align*}
and use \eqref{eq:le3a} for the first term. For the second term, we have
\begin{align*}
    -\sum_j \bg_{j+1/2}^{n+1,\top} \mathcal{L}(\bg_{j+1/2}^{n+1}-\bg_{j+1/2}^n) =&\sum_j (\mathcal{D}^{+}\bg_{j+1/2}^{n+1})^{\top} A^+(\bg_{j+1/2}^{n+1}-\bg_{j+1/2}^n)\\
&+\sum_j (\mathcal{D}^-\bg_{j+1/2}^{n+1})^{\top} A^-(\bg_{j+1/2}^{n+1}-\bg_{j+1/2}^n)\\
=&\sum_j(\bg_{j+1/2}^{n+1}-\bg_{j+1/2}^n)^{\top} (A^+\mathcal{D}^{+}+A^-\mathcal{D}^-) \bg_{j+1/2}^{n+1},
\end{align*}
which proves the lemma.
\end{proof}

\begin{lemma}[Boundedness]\label{le:3.3}
For a given discrete function $\bg_{j+1/2}^n$, the advection operator fulfills the property
    \begin{align*}
        \sum_j [\mathcal{L} \bg_{j+1/2}^{n+1}]^2 \leq 2\sum_j \mathcal{D}^{+}\bg_{j+1/2}^{n+1,\top}\mathbf A^2 \mathcal{D}^{+}\bg_{j+1/2}^{n+1}.
    \end{align*}
\end{lemma}

\begin{proof}
    Note that
    \begin{align*}
        \sum_j [\mathcal{L} &\bg_{j+1/2}^{n+1}]^2\\ =& \sum_j [\mathbf A^+\mathcal{D}^{+} \bg_{j+1/2}^{n+1}]^2 + [\mathbf A^-\mathcal{D}^- \bg_{j+1/2}^{n+1}]^2 + 2 \mathcal{D}^{+} \bg_{j+1/2}^{n+1,\top}\mathbf A^+\mathbf A^- \mathcal{D}^- \bg_{j+1/2}^{n+1}.
    \end{align*}
    For the last term we have with $\mathbf{M}^{\pm} := \frac12\left(\mathbf{M} \pm |\mathbf{M}|\right)$
    \begin{align*}
        \mathcal{D}^{+} \bg_{j+1/2}^{n+1,\top}\mathbf A^+\mathbf A^- \mathcal{D}^- \bg_{j+1/2}^{n+1} = \mathcal{D}^{+} \bg_{j+1/2}^{n+1,\top}\mathbf{T}\mathbf{M}^+\mathbf{T}^{\top}\mathbf{T}\mathbf{M}^-\mathbf{T}^{\top} \mathcal{D}^- \bg_{j+1/2}^{n+1}.
    \end{align*}
    With $\mathbf{T}_{0} := (\sqrt{w_k}p_0(\mu_k))_{k=1}^{N+1} = \left(\sqrt{\frac{w_k}{2}}\right)_{k=1}^{N+1}$ we have $\mathbf{T}^{\top}\mathbf{T} = \mathbf{I}-\mathbf{T}_0\mathbf{T}_0^{\top}$, i.e.,
    \begin{align*}
        \mathcal{D}^{+} \bg_{j+1/2}^{n+1,\top}\mathbf A^+\mathbf A^- \mathcal{D}^- \bg_{j+1/2}^{n+1} = -\mathcal{D}^{+} \bg_{j+1/2}^{n+1,\top}\mathbf{T}\mathbf{M}^+\mathbf{T}_0\mathbf{T}_0^{\top}\mathbf{M}^-\mathbf{T} \mathcal{D}^- \bg_{j+1/2}^{n+1}.
    \end{align*}
    This is a multiplication of two scalars. With $\mathbf{\widehat g}_j:=\mathbf{T}^{\top}\bg_{j+1/2}^{n+1}$ these read
    \begin{align*}
        a_j := \mathcal{D}^{+} \bg_{j+1/2}^{n+1,\top}\mathbf{T}\mathbf{M}^+\mathbf{T}_0 = \mathcal{D}^{+} \mathbf{\widehat g}_j^{\top}\mathbf{M}^+\mathbf{T}_0 = \sum_{k}\mathcal{D}^{+}\widehat g_{jk}\mu_k^+\sqrt{\frac{w_k}{2}},\\
        b_j := -\mathbf{T}_0^{\top}\mathbf{M}^-\mathbf{T} \mathcal{D}^- \bg_{j+1/2}^{n+1} = \mathbf{T}_0^{\top}|\mathbf{M}^-| \mathcal{D}^- \mathbf{\widehat g}_j = \sum_{\ell}\mathcal{D}^-\widehat g_{j\ell}|\mu_{\ell}^-|\sqrt{\frac{w_{\ell}}{2}}.
    \end{align*}
    Hence, for the product of these scalars we have with $a_j b_j \leq \frac12 a_j^2 + \frac12 b_j^2$
    \begin{align*}
        \mathcal{D}^{+} \bg_{j+1/2}^{n+1,\top}\mathbf A^+\mathbf A^- \mathcal{D}^- \bg_{j+1/2}^{n+1} \stackrel{\text{Young}}{\leq} \frac14\sum_{k}(\mathcal{D}^{+}\widehat g_{jk})^2(\mu_k^+)^{2}w_k + \frac14\sum_{\ell}(\mathcal{D}^-\widehat g_{j\ell})^2(\mu_{\ell}^-)^{2}w_{\ell}.
    \end{align*}
    Then, summing over $j$ and using an index shift for the second term gives
    \begin{align*}
         \sum_j \mathcal{D}^{+} \bg_{j+1/2}^{n+1,\top}\mathbf A^+\mathbf A^- \mathcal{D}^- \bg_{j+1/2}^{n+1} \leq & \frac{1}{2}\sum_{k,j}(\mathcal{D}^{+}\widehat g_{jk})^2\mu_k^{2}w_k\\
         =& \frac{1}{2}\sum_j \mathcal{D}^{+}\bg_{j+1/2}^{n+1,\top}\mathbf A^2 \mathcal{D}^{+}\bg_{j+1/2}^{n+1}.
    \end{align*}
    Hence,
    \begin{align*}
        \sum_j [\mathcal{L} \bg_{j+1/2}^{n+1}]^2 \leq 2\sum_j \mathcal{D}^{+}\bg_{j+1/2}^{n+1,\top}\mathbf A^2 \mathcal{D}^{+}\bg_{j+1/2}^{n+1}
    \end{align*}
\end{proof}
The previous Lemmas then allow proving Theorem~\ref{th:estabilityFull}:
\begin{proof}[Proof (Theorem~\ref{th:estabilityFull})]
Multiplying $\rho_j^{n+1}$ to the first equation of (\ref{eq:schemeFull}) and summing over $j$ yields
\begin{equation}  \label{eq:rhorho}
\frac{1}{2}\|\rho^{n+1}\|^2-\frac{1}{2}\|\rho^n\|^2+\frac{1}{2}\|\rho^{n+1}-\rho^n\|^2+\Delta t \Delta x a_0\sum_j \rho_j^{n+1}\mathcal{D}^- g_{1,j+1/2}^{n+1}=0.
\end{equation}
Multiplying $\bg_{j+1/2}^{n+1}$ to the second equation of \eqref{eq:schemeFull} and summing over $j$ yields
\begin{equation}\label{eq:th5-2Eq1}
\begin{split}
&\frac{1}{2}\|\bg^{n+1}\|^2-\frac{1}{2}\|\bg^n\|^2+\frac{1}{2}\|\bg^{n+1}-\bg^n\|^2+\frac{\Delta t \Delta x}{\varepsilon} \sum_j \bg_{j+1/2}^{n+1,\top} \mathcal{L}\bg_{j+1/2}^n\\
&=-\frac{\Delta t}{\varepsilon^2}\sum_j \bg_{j+1/2}^{n+1,\top}\ba (\rho_{j+1}^n-\rho_j^n)-\frac{\Delta t}{\varepsilon^2}\sum_j \sigma_{j+1/2}\bg_{j+1/2}^{n+1,\top}\bg_{j+1/2}^{n+1}\Delta x.
\end{split}
\end{equation}
Using the summation by parts Lemma~\ref{le:sumbyparts} and $\sigma_{j+1/2} \geq \sigma_0 := \min_j \sigma_{j+1/2}$, we obtain
\begin{equation} \label{gg}
\begin{split}
&\frac{1}{2}\|\bg^{n+1}\|^2-\frac{1}{2}\|\bg^n\|^2+\frac{1}{2}\|\bg^{n+1}-\bg^n\|^2+\frac{\Delta t \Delta x}{\varepsilon} \sum_j \bg_{j+1/2}^{n+1,\top} \mathcal{L}\bg_{j+1/2}^n\\
&-\frac{\Delta t \Delta x a_0}{\varepsilon^2} \sum_j \rho_j^n\mathcal{D}^-g_{1,j+1/2}^{n+1}\leq-\frac{\Delta t \sigma_0}{\varepsilon^2}\|\bg^{n+1}\|^2.
\end{split}
\end{equation}
Adding (\ref{eq:rhorho}) and $\varepsilon^2\times$(\ref{gg}), we have 
\begin{equation} \label{19}
\begin{split}
&\frac{1}{2}e^{n+1}-\frac{1}{2}e^n+\frac{1}{2}\|\rho^{n+1}-\rho^n\|^2+\frac{\varepsilon^2}{2}\|\bg^{n+1}-\bg^n\|^2\\
&+\Delta t \Delta x a_0\sum_j ( \rho_j^{n+1}-\rho_j^n)\mathcal{D}^- g_{1,j+1/2}^{n+1}+\varepsilon\Delta t \Delta x\sum_j \bg_{j+1/2}^{n+1,\top} \mathcal{L}\bg_{j+1/2}^n\\
&\leq - \Delta t \sigma_0\|\bg^{n+1}\|^2.
\end{split}
\end{equation}
Due to Lemma~\ref{le:3.2}, equation \eqref{19} becomes
\begin{equation*} 
\begin{split}
&\frac{1}{2}e^{n+1}-\frac{1}{2}e^{n}+\frac{1}{2}\|\rho^{n+1}-\rho^n\|^2+\frac{\varepsilon^2}{2}\|\bg^{n+1}-\bg^n\|^2\\
&+\Delta t \Delta x a_0\sum_j ( \rho_j^{n+1}-\rho_j^n)\mathcal{D}^- g_{1,j+1/2}^{n+1}\\
&+\varepsilon\Delta t \frac{\Delta x^2}{2} \sum_j \mathcal{D}^{+}\bg_{j+1/2}^{n+1,\top}|A|\mathcal{D}^{+}\bg_{j+1/2}^{n+1}+\varepsilon\Delta t \Delta x \sum_j(\bg_{j+1/2}^{n+1}-\bg_{j+1/2}^n)^{\top} \mathcal{L} \bg_{j+1/2}^{n+1}\\
&\leq - \Delta t \sigma_0  \|\bg^{n+1}\|^2.
\end{split}
\end{equation*}
By Young's inequality,
\begin{equation*}
\Delta t \Delta x a_0\sum_j ( \rho_j^{n+1}-\rho_j^n)\mathcal{D}^- g_{1,j+1/2}^{n+1}\leq \frac{1}{2}\|\rho^{n+1}-\rho^n\|^2+\frac{1}{2}\Delta t^2\Delta xa_0^2\sum_j (\mathcal{D}^{+} g^{n+1}_{1,j+1/2})^2
\end{equation*}
as well as
\begin{equation}
\begin{split}
&\varepsilon\Delta t \Delta x \sum_j(\bg_{j+1/2}^{n+1}-\bg_{j+1/2}^n)^{\top} \mathcal{L} \bg_{j+1/2}^{n+1}\leq \frac{\varepsilon^2}{2}\|\bg_{j+1/2}^{n+1}-\bg_{j+1/2}^n\|^2\\
&\quad \quad +\frac{1}{2}\Delta t^2\Delta x \sum_j [\mathcal{L} \bg_{j+1/2}^{n+1}]^2.
\end{split}
\end{equation}
Together with Lemma~\ref{le:3.3} this gives
\begin{equation} 
\begin{split}
&\frac{1}{2}e^{n+1}-\frac{1}{2}e^{n}\\
\leq & \frac{1}{2}\Delta t^2\Delta x a_0^2\sum_j (\mathcal{D}^{+} g^{n+1}_{1,j+1/2})^2+\Delta t^2\Delta x \sum_j \mathcal{D}^{+}\bg_{j+1/2}^{n+1,\top}\mathbf A^2 \mathcal{D}^{+}\bg_{j+1/2}^{n+1}\\
&-\varepsilon\Delta t \frac{\Delta x^2}{2} \sum_j \mathcal{D}^{+}\bg_{j+1/2}^{n+1,\top}|\mathbf A|\mathcal{D}^{+}\bg_{j+1/2}^{n+1}- \Delta t \sigma_0  \|\bg^{n+1}\|^2\\
=&\frac{1}{2}\Delta t\Delta x \sum_j \mathcal{D}^{+}\bg_{j+1/2}^{n+1,\top}\left[\Delta t(2\mathbf A^2+\ba\ba^{\top}) -\varepsilon \Delta x |\mathbf A|\right]\mathcal{D}^{+}\bg_{j+1/2}^{n+1}- \Delta t \sigma_0  \|\bg^{n+1}\|^2.
\end{split}
\end{equation}
Note that with Lemma~\ref{le:AtoAf} and $\mathbf{h}_{j+1/2} := (0,g_{j+1/2,1},\cdots,g_{j+1/2,N})^{\top}$ as well as $\mathbf{\widehat h}_{j+1/2} := \mathbf{T}_f^{\top}\mathbf{h}_{j+1/2}$, we have
\begin{align*}
    \mathcal{D}^{+}\bg_{j+1/2}^{n+1,\top}\mathbf A^2 \mathcal{D}^{+}\bg_{j+1/2}^{n+1} = (\mathcal{D}^{+}\mathbf{\widehat h}_{j+1/2}^{n+1})^{\top}\mathbf M^2 \mathcal{D}^{+}\mathbf{\widehat h}_{j+1/2}^{n+1}\\
    \mathcal{D}^{+}\bg_{j+1/2}^{n+1,\top}|\mathbf A| \mathcal{D}^{+}\bg_{j+1/2}^{n+1} = (\mathcal{D}^{+}\mathbf{\widehat h}_{j+1/2}^{n+1})^{\top}|\mathbf M| \mathcal{D}^{+}\mathbf{\widehat h}_{j+1/2}^{n+1}
\end{align*}
Moreover, according to Lemma~\ref{le:AtoAf}, we have 
\begin{align*}
    \mathcal{D}^{+}\bg_{j+1/2}^{n+1,\top}\ba\ba^{\top} \mathcal{D}^{+}\bg_{j+1/2}^{n+1} = (\mathcal{D}^{+}\mathbf{\widehat h}_{j+1/2}^{n+1})^{\top}\mathbf{T}_f^{\top}\ba_f\ba_f^{\top}\mathbf{T}_f \mathcal{D}^{+}\mathbf{\widehat h}_{j+1/2}^{n+1}.
\end{align*}
Since $\mathbf{T}_f^{\top}\mathbf{a}_f = \frac{1}{\sqrt{3}}\mathbf{T}_f^{\top}\mathbf{e}_2 = \sqrt{\frac{w_k}{3}} p_1(\mu_k) = \sqrt{\frac{w_k}{2}} \mu_k$, we have
\begin{align*}
    (\mathcal{D}^{+}\mathbf{\widehat h}_{j+1/2}^{n+1})^{\top}\mathbf{T}_f^{\top}\ba_f\ba_f^{\top}\mathbf{T}_f \mathcal{D}^{+}\mathbf{\widehat h}_{j+1/2}^{n+1} =& \left(\sum_{k=1}^{N+1} \mathcal{D}^{+}\widehat h_{j+1/2,k}^{n+1}\sqrt{\frac{w_k}{2}} \mu_k\right)^2\\
    =& \frac12\sum_{k,\ell=1}^{N+1} \mathcal{D}^{+}\widehat h_{j+1/2,k}^{n+1}\sqrt{w_k} \mu_k \mathcal{D}^{+}\widehat h_{j+1/2,\ell}^{n+1}\sqrt{w_{\ell}} \mu_{\ell}\\
    \stackrel{\text{Young}}{\leq}& \frac14\sum_{k,\ell=1}^{N+1} \left(\mathcal{D}^{+}\widehat h_{j+1/2,k}^{n+1}\right)^2 w_k \mu_k^2\\
    &+ \frac14\sum_{k,\ell=1}^{N+1} \left(\mathcal{D}^{+}\widehat h_{j+1/2,\ell}^{n+1}\right)^2 w_{\ell} \mu_{\ell}^2\\
    =& \frac{N+1}{2}\sum_{k=1}^{N+1} \left(\mathcal{D}^{+}\widehat h_{j+1/2,k}^{n+1}\right)^2 w_k \mu_k^2.
\end{align*}
Hence, since $\|\bg^{n+1}\|^2 = \|\mathbf{\widehat h}^{n+1}\|^2 = \Delta x\sum_{k,j} \widehat h_{j+1/2,k}^2$ we have that $\Delta e := \frac{1}{2}e^{n+1}-\frac{1}{2}e^{n}$ fulfills 
\begin{align*}
    \Delta e \leq& \frac{\Delta t\Delta x}{2} \sum_{j,k} \left(\mathcal{D}^{+}\widehat h_{j+1/2,k}^{n+1}\right)^2\left[\Delta t(2\mu_k^2+(N+1)w_k \mu_k^2)-\varepsilon \Delta x |\mu_k|\right]\\
    &- \Delta t\Delta x \sigma_0 \sum_{j,k} \left(\widehat h_{j+1/2,k}^{n+1}\right)^2
\end{align*}
Since
\begin{align*}
    \sum_{j} \left(\mathcal{D}^{+}\widehat h_{j+1/2,k}^{n+1}\right)^2 =& \frac{1}{\Delta x^2}\sum_j (\widehat{h}_{j+3/2,k}-\widehat{h}_{j+1/2,k})^2\\
    =& \frac{2}{\Delta x^2}\sum_j \widehat{h}_{j+1/2,k}^2-\frac{2}{\Delta x^2}\sum_j\widehat{h}_{j+1/2,k}\widehat{h}_{j+3/2,k}\\
    \stackrel{\text{Young}}{\leq}& \frac{2}{\Delta x^2}\sum_j \widehat{h}_{j+1/2,k}^2+\frac{1}{\Delta x^2}\sum_j\widehat{h}_{j+1/2,k}+\frac{1}{\Delta x^2}\sum_j\widehat{h}_{j+3/2,k}\\
    =& \frac{4}{\Delta x^2}\sum_j \widehat{h}_{j+1/2,k}^2
\end{align*}
we have
\begin{align*}
    \Delta e \leq 2\frac{\Delta t}{\Delta x} \sum_{j,k} \left(\widehat h_{j+1/2,k}^{n+1}\right)^2\left[\Delta t(2\mu_k^2+(N+1)w_k \mu_k^2) -\varepsilon \Delta x |\mu_k| - \frac12\Delta x^2\sigma_0\right].
\end{align*}
To ensure stability, we must have for all $k$, where $\mu_k\neq 0$
\begin{align*}
    &\Delta t(2\mu_k^2+(N+1)w_k \mu_k^2) \leq \varepsilon \Delta x |\mu_k| + \frac12\Delta x^2\sigma_0\\
    \Leftrightarrow \; & \Delta t(2+(N+1)w_k ) \leq \varepsilon \frac{\Delta x}{|\mu_k|} + \frac{\sigma_0\Delta x^2}{2\mu_k^2}
\end{align*}
Hence, picking 
\begin{align*}
    \Delta t \leq \frac1{2+(N+1)w_k}\left(\varepsilon \frac{\Delta x}{|\mu_k|} + \frac{\sigma_0\Delta x^2}{2\mu_k^2}\right)
\end{align*}
ensures stability. Note that $(N+1)w_k$ remains bounded.
\end{proof}

\section{Dynamical low-rank approximation for modal micro-macro}\label{sec:DLRAmicromacro}
In the following, we derive evolution equations based on the modal micro-macro scheme \eqref{eq:schemeFull}. The DLRA integrator used is the BUG integrator \cite{ceruti2022unconventional}, which preserves the energy stability property of Theorem~\ref{th:estabilityFull} while being asymptotic--preserving.
\subsection{BUG integrator steps}
Let us write down the individual steps for the BUG integrator. Recall the microscopic modal approximation of the P$_N$ system \eqref{eq:PN} which reads 
\begin{align}
    \partial_t \mathbf{g}(t,x) = -\frac{1}{\varepsilon} \mathbf{A} \partial_{x} \mathbf{g}(t,x)  - \frac{1}{\varepsilon^2}\mathbf{a} \partial_{x}\rho - \frac{\sigma}{\varepsilon^2}\mathbf{g}(t,x).
\end{align}
To reduce computational costs and memory requirements, a low-rank ansatz for the microscopic correction is defined as 
\begin{align}
    \mathbf{g}(t,x) \approx \sum_{i,j = 1}^r X_i(t,x)S_{ij}(t) \mathbf{V}_j(t) = \mathbf{X}(t,x)^{\top} \mathbf{S}(t)\mathbf{V}(t)^{\top},
\end{align}
where $\mathbf{S}(t) = (S_{ij}(t))_{i,j = 1}^r\in \mathbb{R}^{r\times r}$, $\mathbf{X}(t,x) = (X_i(t,x))_{i=1}^r\in\mathbb{R}^r$ as well as $\mathbf{V}(t) = (V_{kj}(t))_{k,j=1}^{N,r}\in\mathbb{R}^{N\times r}$.
Then, the substeps of the BUG integrator to integrate the factorized solution from time $t_n$ to $t_{n+1}$ are given as follows:

\textbf{$K$-step}: For an initial condition $\mathbf{K}(t_n,x) = \mathbf{X}^n(x)^{\top} \mathbf{S}^n$, the $L$-step reads
\begin{align}
    \partial_t \mathbf{K}(t,x) = -\frac{1}{\varepsilon} \mathbf{V}^{n,\top} \mathbf{A} \mathbf{V}^n \partial_{x} \mathbf{K}(t,x) - \frac{1}{\varepsilon^2}\mathbf{a}^{\top} \mathbf{V}^n  \partial_{x}\rho - \frac{\sigma}{\varepsilon^2}\mathbf{K}(t,x).
\end{align}
Omitting the dependency on space in our notation, we define the time discretization
\begin{align}\label{eq:KstepDT}
    \mathbf{K}^{n+1} = \mathbf{K}^{n} - \frac{\Delta t}{\varepsilon} \mathbf{V}^{n,\top} \mathbf{A} \mathbf{V}^n \partial_{x} \mathbf{K}^{n} - \frac{\Delta t}{\varepsilon^2}\mathbf{a}^{\top} \mathbf{V}^n \partial_{x}\rho^{n} - \frac{\Delta t\sigma}{\varepsilon^2}\mathbf{K}^{n+1}
\end{align}

\textbf{$L$-step}: For an initial condition $\mathbf{L}(t_n) = \mathbf{V}^n\mathbf{S}^{n,\top}$, the $L$-step reads
\begin{align}
    \mathbf{\dot L}(t) = -\frac{1}{\varepsilon}\mathbf{A}^{\top}\mathbf{L}(t)\langle \partial_{x} \mathbf{X}^n,\mathbf X^{n,\top}\rangle- \frac{1}{\varepsilon^2}\mathbf{a}\langle \partial_{x}\rho,\mathbf X^{n,\top}\rangle - \frac{\sigma}{\varepsilon^2}\mathbf L(t)\langle \mathbf X^n \mathbf{X}^{n,\top}\sigma\rangle.
\end{align}
We define the time discretization
\begin{align}
    \mathbf{L}^{n+1} =\mathbf{L}^{n} -&\frac{\Delta t}{\varepsilon}\mathbf{A}^{\top}\mathbf{L}^n\langle \partial_{x} \mathbf{X}^n,\mathbf X^{n,\top}\rangle \nonumber\\
    -& \frac{\Delta t}{\varepsilon^2}\mathbf{a}\langle \partial_{x}\rho,\mathbf X^{n,\top}\rangle - \frac{\Delta t\sigma}{\varepsilon^2}\mathbf L^{n+1}\langle \mathbf X^n \mathbf{X}^{n,\top}\sigma\rangle.
\end{align}

\textbf{$S$-step}: For an initial condition $\mathbf{S}(t_n) = \mathbf{X}^{n+1,\top}\mathbf{X}^n\mathbf{S}^{n}\mathbf{V}^{n,\top}\mathbf{V}^{n+1}$, the $S$-step reads
\begin{align}
    \mathbf{\dot S} =& -\frac{1}{\varepsilon}\langle \partial_{x} \mathbf{X}^{n+1},\mathbf X^{n+1,\top}\rangle \mathbf{S}(t)  \mathbf{V}^{n+1,\top} \mathbf{A} \mathbf{V}^{n+1} \nonumber\\
    &- \frac{1}{\varepsilon^2}\left(\langle \partial_{x}\rho,\mathbf X^{n+1}\rangle \cdot \left(\mathbf{a}^{\top} \mathbf{V}^{n+1}\right)^{\top}+ \langle \mathbf X^{n+1}\mathbf X^{n+1,\top}\sigma\rangle \mathbf{S}(t)\right).
\end{align}
Defining $\mathbf{\tilde S}^n = \mathbf{X}^{n+1,\top}\mathbf{X}^n\mathbf{S}^{n}\mathbf{V}^{n,\top}\mathbf{V}^{n+1}$, the chosen time discretization is
\begin{align}
\mathbf{S}^{n+1} = \mathbf{\tilde S}^{n}&-\frac{\Delta t}{\varepsilon}\langle \partial_{x} \mathbf{X}^{n+1},\mathbf X^{n+1,\top}\rangle \mathbf{\tilde S}^n  \mathbf{V}^{n+1,\top} \mathbf{A} \mathbf{V}^{n+1} \nonumber\\
    &- \frac{\Delta t}{\varepsilon^2}\left(\langle \partial_{x}\rho^n,\mathbf X^{n+1}\rangle \cdot \left(\mathbf{a}^{\top} \mathbf{V}^{n+1}\right)^{\top}+ \langle \mathbf X^{n+1}\mathbf X^{n+1,\top}\sigma\rangle \mathbf{S}^{n+1}\right).
\end{align}

\subsection{AP-property}
The proposed semi-discrete dynamical low-rank scheme is asymptotic--preserving:
\begin{theorem}
In the limit $\varepsilon\rightarrow 0$, the proposed method preserves the semi-discrete diffusion equation 
\begin{align}
    \frac{1}{\Delta t}\left(\rho^{n+1}-\rho^{n}\right) - \frac13\partial_{x} \left(\frac{1}{\sigma}\partial_{x}\rho^n\right) = 0.
\end{align}
\end{theorem}
\begin{proof}
In the limit, the $K$-step and $L$-step become
\begin{align*}
    \mathbf{a}^{\top} \mathbf{V}^n  \partial_{x}\rho^{n} =& -\sigma \mathbf K^{n+1},\\
    \mathbf{L}^{n+1}\langle \mathbf X^{n}\mathbf X^{n,\top}\sigma\rangle =& -\mathbf{a}\langle \partial_{x}\rho,\mathbf X^{n,\top}\rangle .
\end{align*}
If $\mathbf{L}^{n+1} = \mathbf{V}^{n+1}\mathbf{S}_L$ and $\mathbf{S}_L\langle \mathbf X^{n}\mathbf X^{n,\top}\sigma\rangle$ is invertible, we know that $\frac{1}{\sigma}\partial_{x}\rho^{n}$ and $\mathbf{a}$ lie in the ranges of the time-updated directional and spatial basis sets. For the $S$-step, $\varepsilon\rightarrow 0$ reveals that 
\begin{align}\label{eq:limitS}
    -\langle \partial_{x}\rho^n,\mathbf X^{n+1}\rangle \cdot  \left(\mathbf{a}^{\top} \mathbf{V}^{n+1}\right)^{\top}= \langle\mathbf X^{n+1}\mathbf X^{n+1,\top}\sigma\rangle \mathbf{S}^{n+1}.
\end{align}
Note that with $\frac{1}{\sigma}\partial_{x}\rho^{n} = \langle\frac{1}{\sigma}\partial_{x}\rho^{n},\mathbf X^{n+1,\top}\rangle \mathbf X^{n+1}$ we have
\begin{align*}
    \langle \partial_{x}\rho^n,\mathbf X^{n+1}\rangle = \left\langle \sigma\frac{1}{\sigma}\partial_{x}\rho^{n},\mathbf X^{n+1}\right\rangle = \left\langle \mathbf X^{n+1}\mathbf X^{n+1,\top}\sigma\right\rangle \left\langle \frac{1}{\sigma}\partial_{x}\rho^{n},\mathbf X^{n+1}\right\rangle.
\end{align*}
Hence, the limiting equation \eqref{eq:limitS} becomes
\begin{align*}
    \mathbf{S}^{n+1} = -\left\langle \frac{1}{\sigma}\partial_{x}\rho^{n},\mathbf X^{n+1}\right\rangle \cdot \left(\mathbf{a}^{\top} \mathbf{V}^{n+1}\right)^{\top}
\end{align*}
and since $\frac{1}{\sigma}\partial_{x}\rho^{n}$ and $\mathbf{a}$ lie in the ranges of our basis, scalar multiplication with $\mathbf X^{n+1}$ and $\mathbf V^{n+1}$ directly gives
\begin{align*}
    \mathbf g^{n+1} = -\frac{1}{\sigma}\partial_{x}\rho^{n}\mathbf{a}.
\end{align*}
Plugging this into the density equation proves the theorem.
\end{proof}

\section{Discretization in space}\label{sec:DLRAdisc}
 To discretize the derived evolution equations of the BUG integrator in space, the same strategy as in Section~\ref{sec:discFull} is applied. Throughout this section, we again assume a one-dimensional spatial domain. To provide a better understanding of the properties that are given by the chosen spatial discretization, we first leave time continuous and discuss the fully discrete setting later.
 \subsection{Time continuous scheme}
 First, we discretize the spatial basis as 
 \begin{align*}
     \mathbf{X}_{j+1/2}(t) := \frac{1}{\Delta x} \int_{x_j}^{x_{j+1}}\mathbf{X}(t,x)\, \mathrm{d}x\in\mathbb{R}^r
 \end{align*}
and define $\mathbf{K}_{j+1/2}(t) := \mathbf{X}_{j+1/2}(t)^{\top}\mathbf{S}(t)\in\mathbb{R}^r$. For the time-continuous $K$-step, this gives
\begin{align*} 
\dot{\mathbf{K}}_{j+1/2}(t)+\frac{1}{\varepsilon}\mathcal{L}_K\mathbf{K}_{j+1/2}(t)=-\frac{1}{\varepsilon^2}\ba^{\top}\mathbf{V}^n \mathcal{D}^{+}\rho_j(t)-\frac{\sigma_{j+1/2}}{\varepsilon^2}\mathbf{K}_{j+1/2}(t),
\end{align*}
where we define the streaming operator of the $K$-step as
\begin{align*}
    \mathcal{L}_K\mathbf{K}_{j+1/2} := \mathcal{D}^-\mathbf{K}_{j+1/2}\mathbf{V}^{n,\top}\mathbf A^+\mathbf{V}^n + \mathcal{D}^{+}\mathbf{K}_{j+1/2}\mathbf{V}^{n,\top}\mathbf A^-\mathbf{V}^n.
\end{align*}
For the $L$-step, we obtain
\begin{align*} 
\dot{\mathbf{L}}(t)+\frac{1}{\varepsilon}\mathcal{L}_L\mathbf{L}(t)=&-\frac{1}{\varepsilon^2}\ba \sum_j \mathbf{X}_{j+1/2}^{n,\top}\mathcal{D}^{+}\rho_j(t)-\frac{1}{\varepsilon^2}\sum_j\sigma_{j+1/2}\mathbf{X}_{j+1/2}^{n}\mathbf{X}_{j+1/2}^{n,\top}\mathbf{L}(t),
\end{align*}
where the streaming operator of the $L$-step reads
\begin{align*}
    \mathcal{L}_L\mathbf{L}(t) := \mathbf A^+\mathbf{L}(t)\sum_j\mathcal{D}^-\mathbf{X}^{n}_{j+1/2}\mathbf{X}^{n,\top}_{j+1/2} + \mathbf A^-\mathbf{L}(t)\sum_j\mathcal{D}^{+}\mathbf{X}^n_{j+1/2}\mathbf{X}^{n,\top}_{j+1/2}.
\end{align*}
Lastly, the $S$-step becomes
\begin{align}\label{eq:SstepContTime}
\dot{\mathbf{S}}(t)+\frac{1}{\varepsilon}\mathcal{L}_S\mathbf{S}(t)=&-\frac{1}{\varepsilon^2} \sum_j \mathbf{X}_{j+1/2}^{n+1,\top}\mathcal{D}^{+}\rho_j(t)(\mathbf{a}^{\top}\mathbf{V}^{n+1})^{\top}\nonumber\\
&-\frac{1}{\varepsilon^2}\sum_j \mathbf{X}_{j+1/2}^{n+1}\sigma_{j+1/2}\mathbf{X}_{j+1/2}^{n+1,\top}\mathbf{S}(t),
\end{align}
where with $\mathbf{\widehat D}^{\pm} := \sum_j\mathbf{X}^{n+1}_{j+1/2} D^{\pm}\mathbf{X}_{j+1/2}^{n+1,\top}\in\mathbb{R}^{r\times r}$ the streaming operator of the $S$-step reads
\begin{align*}
    \mathcal{L}_S\mathbf{S} := \mathbf{\widehat D}^{-}\mathbf{S}\mathbf{V}^{n+1,\top}\mathbf A^+\mathbf{V}^{n+1} + \mathbf{\widehat D}^{+}\mathbf{S}\mathbf{V}^{n+1,\top}\mathbf A^-\mathbf{V}^{n+1}.
\end{align*}
Combined with the macro equation
\begin{align*}
    \dot \rho_j(t)+a_0\mathcal{D}^-\mathbf{X}_{j+1/2}^{n+1,\top}\mathbf{S}^{n+1}\mathbf{V}^{n+1,\top}\mathbf{e}_1=0,
\end{align*}
the derived time-continuous scheme dissipates the energy, as we show in the following theorem.
\begin{theorem}\label{th:eStableSemiDisc}
    The proposed time-continuous method is energy stable. More precisely, for the Frobenius norm $\Vert\cdot\Vert$ we have
        \begin{align*}
        \frac12\frac{d}{dt}\left(\Vert\bm{\rho}\Vert^2 + \varepsilon^2\Vert\mathbf{g}^S\Vert^2 \right)  \leq -\frac{\sigma_{0}(x_{j+1/2})}{\varepsilon^2} \Vert\mathbf{g}^S\Vert^2,
    \end{align*}
    where $\mathbf{g}^S(t):= \mathbf{X}^{n+1}\mathbf{S}(t)\mathbf{V}^{n+1,\top}$ is the solution of the Galerkin step.
\end{theorem}
\begin{proof}
    First, let us multiply the $S$-step \eqref{eq:SstepContTime} with $\mathbf{X}^{1,\top}$ from the left and $\mathbf{V}^1$ from the right, define $\mathbf{g}^S(t):= \mathbf{X}_1\mathbf{S}(t)\mathbf{V}_1^{\top}$ as well as $\mathbf{P}^X\in\mathbb{R}^{N_x \times N_x}$ with entries $P_{\alpha j}^X:= \sum_k X_{\alpha+1/2, k}^{n+1}X_{j+1/2, k}^{n+1}$ and $\mathbf{P}^V:= \mathbf{V}^{n+1}\mathbf{V}^{n+1,\top}\in\mathbb{R}^{N \times N}$. Then, in index notation when using Einstein's sum convention we have
    \begin{align*}
    \dot{g}^S_{\alpha +1/2, \beta}(t) =& -\frac{1}{\varepsilon}P^X_{\alpha j}\mathcal{D}^-g^S_{j+1/2,\ell} A^+_{k\ell}P^V_{k\beta}-\frac{1}{\varepsilon}P^X_{\alpha j}\mathcal{D}^{+}g^S_{j+1/2,\ell}A^{-}_{k\ell}P^V_{k\beta}\\
    &-\frac{1}{\varepsilon^2} P^X_{\alpha j}\sigma_{s,j+1/2} g^S_{j+1/2,k}P^V_{k\beta} - \frac{1}{\varepsilon^2}P^X_{\alpha j}\mathcal{D}^{+}\rho_ja_kP^V_{k\beta}.
    \end{align*}
    Now we multiply with $g^S_{\alpha+1/2, \beta}$ and sum over $\alpha$ and $\beta$. For this, note that 
    \begin{align}
        P^X_{\alpha j}g^S_{\alpha+1/2, k} = g^S_{j+1/2, k} \quad \text{and} \quad P^V_{k\beta} g^S_{j+1/2, \beta} = g^S_{j+1/2, k}.
    \end{align}
    Moreover, we know from Lemma~\ref{le:3.2} that
    \begin{align*}
        -\sum \mathbf{g}_{j+1/2}^{S,\top}\mathcal{L}\mathbf{g}^S_{j+1/2} =& -g^S_{j+1/2,k}\mathcal{D}^-g^S_{j+1/2,\ell} A^+_{k\ell}- g^S_{j+1/2,k}\mathcal{D}^{+}g^S_{j+1/2,\ell}A^{-}_{k\ell}\\
        =& -\mathcal{D}^{+}g_{j+1/2,k}^{n+1}| A|_{k\ell}\mathcal{D}^{+}\bg_{j+1/2,\ell}^{n+1} = -\left(\mathcal{D}^{+}g_{j+1/2,k}^{n+1}| A|_{k\ell}^{1/2}\right)^2.
    \end{align*}
    Hence, we directly have
    \begin{align*}
        \frac12\frac{d}{dt}\Vert\mathbf{g}^S\Vert^2 =& -\frac1{\varepsilon}\left(\mathcal{D}^{+}g_{j+1/2,k}^{n+1}| A|_{k\ell}^{1/2}\right)^2-\frac{\sigma_{s,j+1/2}}{\varepsilon^2} \left(g_{j+1/2,k}^S\right)^2 - \frac{1}{\varepsilon^2}g_{j+1/2,k}^S\mathcal{D}^{+}\rho_j a_k.
    \end{align*}
    Then, we multiply the macroscopic equation with $\rho_{j}$ and sum over $j$ which yields
    \begin{align*}
        \frac12\frac{d}{dt}\Vert\bm{\rho}\Vert^2 = -\rho_{j} \mathcal{D}^{+}g^S_{j+1/2,k}a_{k} 
        =g_{j+1/2,k}^S a_{k}\mathcal{D}^{+}\rho_{j}.
    \end{align*}
    All together, we then have
    \begin{align*}
       \frac12\frac{d}{dt}\left(\Vert\bm{\rho}\Vert^2 + \varepsilon^2\Vert\mathbf{g}^S\Vert^2 \right) = -\frac1{\varepsilon}\left(\mathcal{D}^{+}g_{j+1/2,k}^{n+1}| A|_{k\ell}^{1/2}\right)^2 -\frac{\sigma_{s,j+1/2}}{\varepsilon^2} \left(g_{j+1/2,k}^S\right)^2
    \end{align*}
\end{proof}

\subsection{Fully discrete scheme}
Let us now derive a time discretization which recovers a discrete counterpart of Theorem~\ref{th:eStableSemiDisc}. For this, we choose an IMEX scheme to treat scattering terms implicitly. Then, the $K$-steps reads
\begin{align*} 
\mathbf{K}_{j+1/2}^{n+1} = \mathbf{K}_{j+1/2}^{n} -\frac{\Delta t}{\varepsilon}\mathcal{L}_K\mathbf{K}_{j+1/2}^{n}-\frac{\Delta t}{\varepsilon^2}\ba^{\top}\mathbf{V}^n \mathcal{D}^{+}\rho_j^{n}-\frac{\Delta t\sigma_{j+1/2}}{\varepsilon^2}\mathbf{K}_{j+1/2}^{n+1}.
\end{align*}
For the $L$-step, we obtain
\begin{align*} 
\mathbf{L}^{n+1} = \mathbf{L}^{n}-\frac{\Delta t}{\varepsilon}\mathcal{L}_L\mathbf{L}^n-\frac{\Delta t}{\varepsilon^2}\ba \sum_j \mathbf{X}_{j+1/2}^{n,\top}\mathcal{D}^{+}\rho_j^n-\frac{\Delta t}{\varepsilon^2}\sum_j\sigma_{j+1/2}\mathbf{X}_{j+1/2}^{n}\mathbf{X}_{j+1/2}^{n,\top}\mathbf{L}^{n+1},
\end{align*}
The $S$-step with $\mathbf{\widetilde S}^{n} := \sum_j \mathbf X_{j+1/2}^n\mathbf X_{j+1/2}^{n+1}\mathbf{S}^{n}\mathbf{V}^{n,\top}\mathbf{V}^{n+1}$ becomes
\begin{equation}\label{eq:SstepDiscrTime}
\begin{split}
\mathbf{S}^{n+1}=\mathbf{\widetilde S}^{n} &-\frac{\Delta t}{\varepsilon}\mathcal{L}_S\mathbf{\widetilde S}^{n} -\frac{\Delta t}{\varepsilon^2} \sum_j \mathbf{X}_{j+1/2}^{n+1,\top}\mathcal{D}^{+}\rho_j^{n}(\mathbf{a}^{\top}\mathbf{V}^{n+1})^{\top}\\
&-\frac{\Delta t}{\varepsilon^2}\sum_j \mathbf{X}_{j+1/2}^{n+1}\sigma_{j+1/2}\mathbf{X}_{j+1/2}^{n+1,\top}\mathbf{S}^{n+1}.
\end{split}
\end{equation}
Lastly, the macroscopic time update reads
\begin{align}\label{eq:macroDLRA}
    \rho_j^{n+1} = \rho_j^{n} - \Delta ta_0\mathcal{D}^-\mathbf{X}_{j+1/2}^{n+1,\top}\mathbf{S}^{n+1}\mathbf{V}^{n+1,\top}\mathbf{e}_1.
\end{align}
It turns out that the use of the unconventional integrator preserves energy stability under the same CFL condition that we need to show the corresponding result for the full problem.
\begin{theorem}[Energy stability, DLRA]\label{th:estabilityFullDLRA} 
Assume that the time step size $\Delta t$ fulfills the CFL condition \eqref{eq:CFLcondition} of the full scheme, that is
\begin{align*}
    \Delta t \leq C\left(\varepsilon \Delta x + \sigma_0\Delta x^2\right).
\end{align*}
Then, the fully discrete DLRA scheme is energy stable, i.e., $e^{n+1}\leq e^n$.
\end{theorem}
\begin{proof}
    We follow the proof of Theorem~\ref{th:eStableSemiDisc} and multiply the $S$-step \eqref{eq:SstepDiscrTime} with $\mathbf{X}^{1,\top}$ from the left and $\mathbf{V}^1$ from the right. Let us define  $\mathbf{\widetilde g}^{n}:= \mathbf{X}^{n+1}\mathbf{\widetilde S}^n\mathbf{V}^{n+1,\top}$ and $\mathbf{\widetilde g}^{n+1}:= \mathbf{X}^{n+1}\mathbf{S}^{n+1}\mathbf{V}^{n+1,\top}$. Then, in index notation when using Einstein's sum convention we have
    \begin{align*}
     \widetilde{g}_{\alpha +1/2, \beta}^{n+1} = \widetilde{g}_{\alpha +1/2, \beta}^{n} &-\frac{\Delta t}{\varepsilon}P^X_{\alpha j}\mathcal{D}^-\widetilde{g}^{n}_{j+1/2,\ell} A^+_{k\ell}P^V_{k\beta}-\frac{\Delta t}{\varepsilon}P^X_{\alpha j}\mathcal{D}^{+}\widetilde{g}^{n}_{j+1/2,\ell}A^{-}_{k\ell}P^V_{k\beta}\\
    &-\frac{\Delta t}{\varepsilon^2} P^X_{\alpha j}\sigma_{s,j+1/2} \widetilde{g}^{n+1}_{j+1/2,k}P^V_{k\beta} - \frac{\Delta t}{\varepsilon^2}P^X_{\alpha j}\mathcal{D}^{+}\rho_j^na_kP^V_{k\beta}.
    \end{align*}
    Now we multiply with $\widetilde g^{n+1}_{\alpha+1/2, \beta}$ and sum over $\alpha$ and $\beta$. Again note that
    \begin{align}
        P^X_{\alpha j}\widetilde g^{n+1}_{\alpha+1/2, k} = \widetilde g^{n+1}_{j+1/2, k} \quad \text{and} \quad P^V_{k\beta} \widetilde g^{n+1}_{j+1/2, \beta} = \widetilde g^{n+1}_{j+1/2, k}.
    \end{align}
    Then, we obtain
    \begin{equation}\label{eq:th5-2Eq1_DLRA}
        \begin{split}
        &\frac{1}{2}\|\mathbf{\widetilde g}^{n+1}\|^2-\frac{1}{2}\|\mathbf{\widetilde g}^n\|^2+\frac{1}{2}\|\mathbf{\widetilde g}^{n+1}-\mathbf{\widetilde g}^n\|^2+\frac{\Delta t \Delta x}{\varepsilon} \sum_j \mathbf{\widetilde g}_{j+1/2}^{n+1,\top} \mathcal{L}\mathbf{\widetilde g}_{j+1/2}^n\\
        &=-\frac{\Delta t}{\varepsilon^2}\sum_j \mathbf{\widetilde g}_{j+1/2}^{n+1,\top}\ba (\rho_{j+1}^n-\rho_j^n)-\frac{\Delta t}{\varepsilon^2}\sum_j \sigma_{j+1/2}\mathbf{\widetilde g}_{j+1/2}^{n+1,\top}\mathbf{\widetilde g}_{j+1/2}^{n+1}\Delta x.
        \end{split}
    \end{equation}
    Note that this equation is equal to \eqref{eq:th5-2Eq1} when replacing $\mathbf g^n$ with $\mathbf{\widetilde g}^{n}$. Similarly, multiplying $\rho_j^{n+1}$ to the macroscopic equation \eqref{eq:macroDLRA} and summing over $j$ yields
\begin{equation}\label{eq:rhorhoDLRA}
\frac{1}{2}\|\rho^{n+1}\|^2-\frac{1}{2}\|\rho^n\|^2+\frac{1}{2}\|\rho^{n+1}-\rho^n\|^2+\Delta t \Delta x a_0\sum_j \rho_j^{n+1}\mathcal{D}^- \widetilde g_{1,j+1/2}^{n+1}=0,
\end{equation}
which is equivalent to \eqref{eq:rhorho} when again replacing $\mathbf g^{n+1}$ with $\mathbf{\widetilde g}^{n+1}$. Hence, adding \eqref{eq:rhorhoDLRA} and $\varepsilon^2\times$\eqref{eq:th5-2Eq1_DLRA}, we obtain equation \eqref{19}. The remainder of the proof follows the proof of Theorem \ref{th:estabilityFull}.
\end{proof}
\section{Numerical results}\label{sec:numres}
\begin{figure}[h!]
\centering
	\begin{subfigure}{0.49\linewidth}
		\centering
		\includegraphics[width=0.99\linewidth]{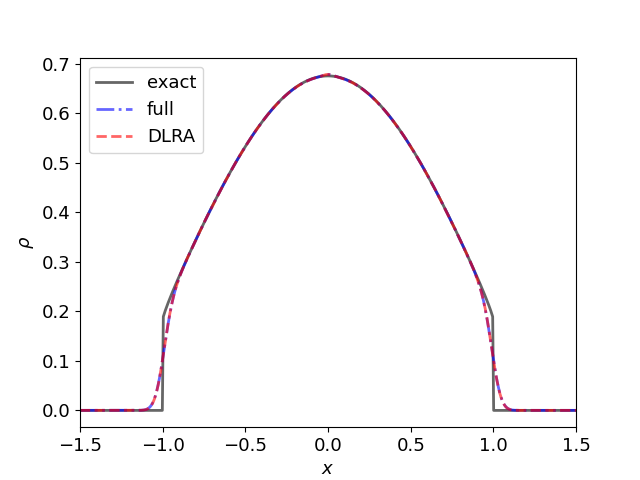}
		\caption{scalar flux, $\varepsilon = 1$}
		\label{fig:rhoEps1PS}
	\end{subfigure}
	\begin{subfigure}{0.49\linewidth}
		\centering
		\includegraphics[width=0.99\linewidth]{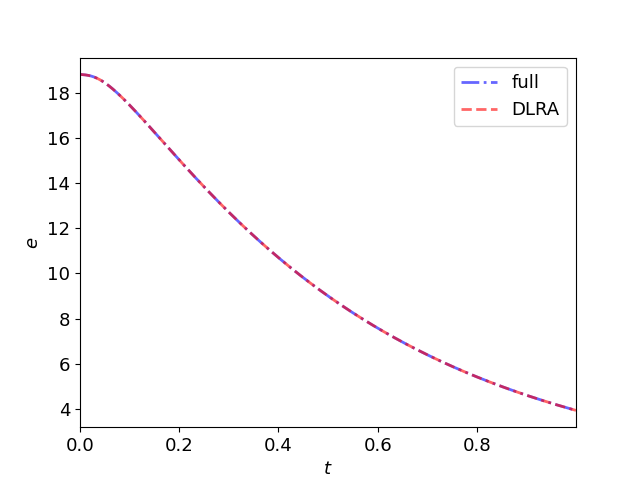}
		\caption{energy dissipation, $\varepsilon = 1$}
		\label{fig:eEps1PS}
	\end{subfigure}
		\begin{subfigure}{0.49\linewidth}
		\centering
		\includegraphics[width=0.99\linewidth]{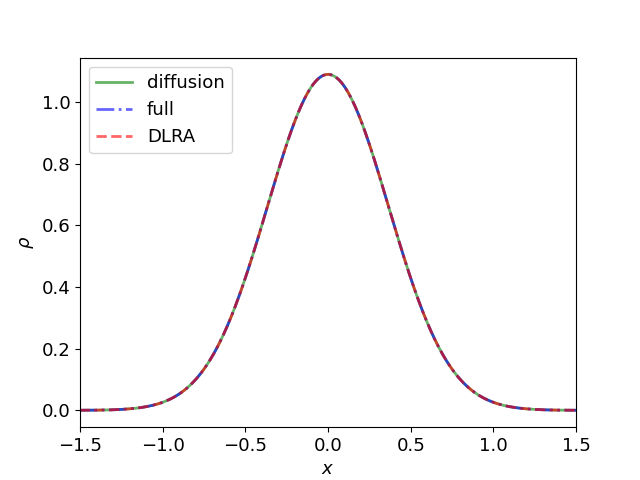}
		\caption{scalar flux, $\varepsilon = 10^{-5}$}
		\label{fig:rhoEps1e-5PS}
	\end{subfigure}
	\begin{subfigure}{0.49\linewidth}
		\centering
		\includegraphics[width=0.99\linewidth]{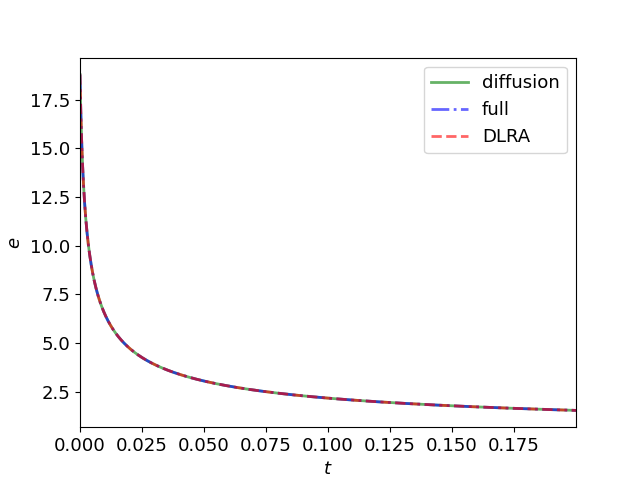}
		\caption{energy dissipation, $\varepsilon = 10^{-5}$}
		\label{fig:eEps1e-5PS}
	\end{subfigure}
	\caption{Scalar flux and energy dissipation.}
	\label{fig:Figure2}
\end{figure}

The following numerical results can be reproduced with the openly available source code \cite{code}. We test the proposed scheme on the plane source testcase~\cite{ganapol2008analytical} in diffusive scaling. For this, we investigate the one-dimensional radiation transport equation in slab geometry in the spatial domain $D = [-1.5,1.5]$ using the initial condition
\begin{align}
    f(t_0) = \frac{1}{\sqrt{2\pi}\sigma}\exp\Big(-\frac{x^2}{2\sigma^2} \Big).
\end{align}
The chosen initial Gaussian has a standard deviation $\sigma = 3\cdot10^{-2}$. That is, the testcase considers particles which are initially positioned around $x=0$ and have an isotropic velocity distribution. As time progresses, particles stream into all directions while undergoing isotropic collisions at a rate of $\sigma = 1$. We use $N=100$ moments to represent the microscopic solution as well as $N_x = 502$ spatial cells. For a Knudsen number of $\varepsilon = 1$, we choose a rank of $r=20$ as well as a final time $t_{\mathrm{end}}=1$. In this case, an analytic solution can be computed according to \cite{ganapol2008analytical}. The chosen time step restriction follows the CFL condition \eqref{eq:CFLcondition} according to
\begin{align*}
    \Delta t = \min_k\left\{\frac1{2+(N+1)w_k}\left(\frac{\varepsilon\Delta x}{|\mu_k|} + \frac{\sigma_0\Delta x^2}{2\mu_k^2}\right)\right\}.
\end{align*}
where for $N=100$ the time step is minimal for $w_k \approx 0.01776$ and $\mu_k \approx -0.81890$.
The resulting macroscopic scalar flux can be found in Figure~\ref{fig:rhoEps1PS}. It is observed that the DLRA solution agrees well with the full P$_N$ solution. Taking a look at the energy dissipation in Figure~\ref{fig:eEps1PS}, we see that both the full as well as the DLRA method exhibit the same energy dissipation. 
For a Knudsen number of $\varepsilon = 10^{-5}$, we pick a final time of $t_{\mathrm{end}} = 0.2$ as well as a rank of $r=3$. In this setting, the time step size is minimal for $w_k \approx 0.01278$ and $\mu_k \approx -0.91079$. Figure~\ref{fig:rhoEps1e-5PS} shows the resulting scalar flux and the corresponding energy dissipation is depicted in Figure~\ref{fig:eEps1e-5PS}. Again, DLRA agrees well with the full rank solution why the method dissipates energy.
\section{Conclusion}
In this work, we derived an asymptotic--preserving dynamical low-rank method which is energy stable under a time step restriction which captures the hyperbolic and parabolic regimes of the radiation transport equation. The proof of energy stability uses the special choice of the stabilization as well as the properties of the ``unconventional'' basis update \& Galerkin step integrator.
\section*{Acknowledgments}
The work of J.\ Kusch was funded by the Deutsche Forschungsgemeinschaft (DFG, German Research Foundation) – 491976834. The work of J.\ Hu was partially supported by NSF DMS-2153208, AFOSR FA9550-21-1-0358 and DOE DE-SC0023164.

\bibliographystyle{siamplain}
\bibliography{references,hu_bibtex}
\end{document}

%% file: art_shared.tex
\usepackage{lipsum}
\usepackage{amsfonts}
\usepackage{graphicx}
\usepackage{epstopdf}
\usepackage{algorithmic}
\ifpdf
  \DeclareGraphicsExtensions{.eps,.pdf,.png,.jpg}
\else
  \DeclareGraphicsExtensions{.eps}
\fi

\newsiamremark{remark}{Remark}
\newsiamremark{hypothesis}{Hypothesis}
\crefname{hypothesis}{Hypothesis}{Hypotheses}
\newsiamthm{claim}{Claim}

\headers{Asymptotic--preserving and energy stable DLRA}{}

\title{Asymptotic--preserving and energy stable dynamical low-rank approximation}

\author{Lukas Einkemmer\thanks{University of Innsbruck,
  (\email{lukas.einkemmer@uibk.ac.at}).}
\and Jingwei Hu\thanks{Department of Applied Mathematics, University of Washington, USA
  (\email{hujw@uw.edu}).}
\and Jonas Kusch\thanks{University of Innsbruck,
  (\email{jonas.kusch1@gmail.com}})}

\usepackage{amsopn}
